\documentclass[12pt]{amsart}

\usepackage{amssymb}

\usepackage{enumerate}

\usepackage{graphicx}

\makeatletter
\@namedef{subjclassname@2020}{%
  \textup{2020} Mathematics Subject Classification}
\makeatother

\usepackage[T1]{fontenc}

\newtheorem{theorem}{Theorem}[section]

\newtheorem{lemma}[theorem]{Lemma}
\newtheorem{proposition}[theorem]{Proposition}

\theoremstyle{definition}
\newtheorem{definition}[theorem]{Definition}

\newtheorem{example}[theorem]{Example}

\numberwithin{equation}{section}

\def\M{\mathcal{M}}
\def\esup{\textup{ess sup}}
\newcommand{\N}{{\mathbb N}}
\newcommand{\al}{\alpha}
\def\dist{\textup{dist}}
\def\d{\textup{d}}
\def\diam{\textup{diam}}
\def\e{\textup{e}}

\begin{document}


\baselineskip=17pt


\title[Almost-compact and compact embedding \dots]{ Almost-compact and compact embeddings of variable exponent spaces}

\author[D.E.Edmunds]{D. E. Edmunds}
\address{Department of mathematics, Pevensey 2 Building, University of Sussex, Brighton
BN1 9QH, UK}
\email{davideedmunds@aol.com}

\author[A.Gogatishvili]{A. Gogatishvili}
\address{Institute of Mathematics CAS \\
\v Zitn\'a 25\\115 67 Praha 1\\Czech Republic}
\email{gogatish@math.cas.cz}

\author[A.Nekvinda]{A. Nekvinda}
\address{Department of Mathematics\\ Faculty of Civil Engineereng\\
       Czech Technical University\\Th\' akurova 7\\16629 Prague 6\\
       Czech Republic}
\email{nekvinda@fsv.cvut.cz}

\date{}

\begin{abstract}
Let $\Omega $ be an open subset of $\mathbb{R}^{N}$, and let $p,\, q:\Omega
\rightarrow \left[ 1,\infty \right] $ be measurable functions. We give a
necessary and sufficient condition for the embedding of the variable
exponent space $L^{p(\cdot )}\left( \Omega \right) $ in $L^{q(\cdot )}\left(
\Omega \right) $ to be almost compact. This leads to a condition on $\Omega, \, p$ and $q$ sufficient to ensure that the Sobolev space $W^{1,p(\cdot
)}\left( \Omega \right) $ based on $L^{p(\cdot )}\left( \Omega \right) $ is
compactly embedded in $L^{q(\cdot )}\left( \Omega \right) ;$ compact
embedding results of this type already in the literature are included as
special cases.
\end{abstract}



\maketitle

\section{Introduction}

Let $\Omega $ be an open subset of $\mathbb{R}^{N}$ and consider the Lebesgue measure on $\Omega$.
If $M\subset \Omega$ is measurable, we write $|M|$ for its measure. Let $p,\, q:\Omega
\rightarrow \left[ 1,\infty \right] $ be measurable. Much attention has been
paid in recent years to the variable exponent space $L^{p(\cdot )}\left(
\Omega \right) ,$ the space $W^{1,p(\cdot )}( \Omega ) $ of
Sobolev type based on $L^{p(\cdot )}( \Omega ) $ and conditions
under which $W^{1,p(\cdot )}( \Omega ) $ is embedded in $%
L^{q(\cdot )}( \Omega ) :$ we refer to \cite{DHHR,ER,ER2} for a comprehensive
account of such matters. The compactness of such an embedding is addressed
here (see  also in the resent paper \cite{FGNR}): we give conditions that are sufficient to ensure compactness yet weak
enough for much earlier work on this topic to be included. To do this we
first establish necessary and sufficient conditions for the embedding of $%
L^{p(\cdot )}( \Omega ) $ in $L^{q(\cdot )}(\Omega) $
to be almost compact.

Let $\mathcal{M}( \Omega ) $ be the family of all measurable functions $u:\Omega
\rightarrow [ -\infty ,\infty ]$; denote by $\chi _{E}$ the
characteristic function of a set $E\subset \Omega $. Given any sequence $
\{ E_{n}\} $ of measurable subsets of $\Omega,$ we write $
E_{n}\rightarrow \emptyset $  a.e. if the characteristic  functions $\chi _{E_{n}}$ converge to $0$ pointwise  almost
everywhere in $\Omega$.
We recall the definition of a Banach function
space: see, for example, \cite{BS}. A normed linear space
$(X,\,\|.\|_X)$ is a Banach function space (BFS for short)
if the following
conditions are satisfied:
\begin{align}
&\text{the norm }\|u\|_X\text{ is defined for all }
            u\in \M (\Omega),\label{2.1}\\
&\text{ and } u\in X\text{ if and only if }\text{} \|u\|_X<\infty;\nonumber\\
&\|u\|_X=\| \ |u| \ \|_X\text{ for every }u\in \M
(\Omega);\label{2.2}\\
&\text{ if }0\le u_n \nearrow u\text{ a.e. in }\Omega, \text{ then } \|u_n\|_X\nearrow \|u\|_X;\label{2.3}\\
&\text{ if }E \subset \Omega \text{ is a measurable set of finite measure, then }\chi_E
              \in X;\label{2.4}\\
&\text{ for every measurable set }E\subset \Omega \text{ of finite measure }|E|
\text{, there exists}\label{2.5}\\
&\text{a positive constant }C_E \text{ such that }\int_E |u(x)| dx \le C_E \|u\|_X\nonumber\\
&\text{ for all }u\in X.\nonumber
\end{align}
If $X$ and $Y$ are Banach function
spaces, then $X$ is said to be almost-compactly embedded in $Y$ and we write $X\overset{\ast}{\hookrightarrow} Y$ if, for every
sequence $\left( E_{n}\right) _{n\in \mathbb{N}}$ of measurable
subsets of $\Omega$ such that $E_{n}\rightarrow \emptyset $  a.e., we
have
\begin{equation*}
\lim_{n\rightarrow \infty }\sup_{\left\Vert u\right\Vert _{X}\leq
1}\left\Vert u\chi _{E_{n}}\right\Vert _{Y}=0.
\end{equation*}%
We believe this notion to have independent interest. Moreover, as we
know from \cite{S}, almost compactness results quickly lead to assertions
concerning the compactness of Sobolev embeddings.

To explain in a little more detail what is achieved,
let $\Omega $
be bounded, $p\in C\left( \overline{\Omega }\right) $ and suppose  that for all $x\in
\Omega, \,\,1<p_{-}\leq p(x)\leq p_{+}<N$ and set
\begin{align}
&p^{\#}(x)=\frac{Np(x)}{N-p(x)};\label{jvhbjovhobvh}
\end{align}
denote by $W_0^{1,p(\cdot )}\left( \Omega
\right) $ the closure of $C_{0}^{\infty }( \Omega ) $ in $%
W^{1,p(\cdot )}( \Omega ) .$ Let $I_{p,q}$ (resp. $I_{p,q,0})$
stand for the embedding of $W^{1,p(\cdot )}( \Omega ) $ $\left(
\text{resp.}\quad W_0^{1,p(\cdot )}( \Omega ) \right) $ in
$L^{q(\cdot )}( \Omega ) .$ Then it is known (see \cite{KR}) that $I_{p,q,0}$ is
compact if there exists $\varepsilon >0$ such that $q(x)\leq p^{\sharp
}(x)-\varepsilon $ for all $x\in \Omega .$ In \cite{KS} the compactness of $%
I_{2,q,0}$ is studied under more general assumptions: it is supposed that
there exists $x_{0}\in \Omega $, a small $\eta >0$, $l\in(0,1)$ and $C>0$  such that $q\left( x_{0}\right) =2N/(N-2)$
and
\begin{equation*}
q(x)\leq \frac{2N}{N-2}-\frac{C}{\left( \log \frac{1}{\left\vert
x-x_{0}\right\vert }\right) ^{l}}
\end{equation*}%
holds for a.e. $x\in \Omega$ with $|x-x_0|\le \eta$. Some generalizations of these results are given in \cite{JPY} and \cite{MOSS}. In
\cite{MOSS} it is assumed that $q(x)=p^{\sharp }(x)$ on a compact set $K,$ and
compactness of $I_{p,q,0}$ is established under some restrictions on $K$ and
on the behavior of $p^{\sharp }(x)-q(x)$ far from $K.$ The principal aim of
the present paper is to establish compactness of $I_{p,q}$ for a wider class of
sets $K$ on which $q$ is allowed to have the same values as $p^{\sharp }:$ various
examples of Cantor type are given for which this is possible.

First we find a necessary and sufficient condition for the embedding of $L^{p(\cdot)}(\Omega)$ into $L^{q(\cdot)}(\Omega)$ to be almost compact. As an application we establish the compactness of the Sobolev embedding mentioned above under more general conditions than those previously available.

\section{Preliminaries}

Let $X$ and $Y$ be Banach function spaces on a bounded open set $\Omega$ of $\mathbb{R}^N$ with norms $\left\Vert \cdot \right\Vert
_{X},\left\Vert \cdot \right\Vert _{Y}$ respectively. We say that $X$ is
embedded in $Y,$ and write $X\hookrightarrow Y,$ if there exists $c>0$ such
that $\left\Vert u\right\Vert _{Y}\leq c\left\Vert u\right\Vert _{X}$ for
all $u\in X.$ The space $X$ is said to be compactly embedded in $Y,$ and we
write $X\hookrightarrow \hookrightarrow Y,$ if given any sequence $\left\{
u_{n}\right\} _{n\in \mathbb{N}}$ such that $\left\Vert u_{n}\right\Vert
_{X}\leq 1$ for all $n\in\mathbb{N}$, there are $\left\{ u_{n(k)}\right\}_{k\in \mathbb{N}} \subset Y$ and
a point $u\in Y$ such that $\left\Vert u_{n(k)}-u\right\Vert _{Y}\rightarrow
0.$

%
\begin{definition}
Let $X$ be a BFS. The Sobolev space $W^1(X)$ is defined to be the set of all weakly differentiable functions $u\in\M(\Omega)$ with
\begin{align*}
&\|u\|_{W^1(X)}=\|u\|_X+\| \nabla u \|_X<\infty.
\end{align*}
\end{definition}

The following proposition is proved in \cite{S}, see Theorem 3.2.

\begin{proposition}\label{wiuhvuvhbhobi}
Let $X,Y,Z$ be BFSs and assume
\begin{align*}
&W^1(X)\hookrightarrow Y,\ \ Y\overset{\ast}{\hookrightarrow} Z.
\end{align*}
Then
\begin{align*}
&W^1(X)\hookrightarrow\hookrightarrow Z.
\end{align*}
\end{proposition}

We now define variable Lebesgue spaces. Let $\mathcal{E}(\Omega)$ denote the set of all measurable functions $p(\cdot):\Omega\rightarrow[1,\infty]$. Given $p(\cdot)\in \mathcal{E}(\Omega)$, we set $\Omega_1 =\{x\in\Omega;\quad p(x)=1\}$,  $\Omega_\infty=\{x\in\Omega;\quad p(x)=\infty\}$ and  $\Omega_0=\Omega \setminus(\Omega_1\cup\Omega_\infty)$.
We adopt the notation
\begin{align*}
&p_{-}=\inf\{p(x); \,\, x\in\Omega_0\},\, \, \, p_{+}=\sup\{p(x); \,\, x\in\Omega_0\}\quad \text{if}\quad  |\Omega_0|>0
\end{align*}
and $p_{-}=p_{+}=1$ if  $|\Omega_0|=0$,
\begin{align*}
&p'(x)=
\begin{cases}
    \infty & \text{ for }p(\cdot)\in\Omega_1,\\
    1 &  \text{ for }p(\cdot)\in \Omega_\infty,\\
    \frac{p(x)}{p(x)-1} &   \text{ for } p(\cdot)\in \Omega_0,
\end{cases}
\end{align*}
\begin{equation}\label{cr} c_{p(\cdot)}=\|\chi_{{}_{\Omega_1}}\|_\infty+ \|\chi_{{}_{\Omega_0}}\|_\infty+\|\chi_{{}_{\Omega_\infty}}\|_\infty,  \quad r_{p(\cdot)}=
c_{p(\cdot)}+\frac{1}{p_{-}}- \frac{1}{p_{+}}.
\end{equation}

For a function $u:\Omega\rightarrow\mathbb{R}$ we define a modular
\begin{align}
&m_{p(\cdot)}(u)=\int_{\Omega\setminus \Omega_\infty} |u(x)|^{p(x)}dx+\esup_{x\in\Omega_\infty}|u(x)|\label{efjkvefivhbobh}.
\end{align}
The space $L^{p(\cdot)}(\Omega)$ is then defined as the set of all measurable functions $u$ on $\Omega$ with a finite norm
$$
\|u\|_{p(\cdot)}:=\inf\Big\{\lambda>0;\,\, m_{p(\cdot)}(u/\lambda)\leq 1\Big\}.
$$

\begin{theorem}[H\"older's inequality, see \cite{KR}, Theorem 2.1]
Let $p(\cdot)\in \mathcal{E}(\Omega)$. Then the inequality
\begin{align*}
&\int_\Omega |f(x)g(x)|dx\le r_{p(\cdot)}\|f\|_{p(\cdot)}\|g\|_{p'(\cdot)},
\end{align*}
 holds for every $f\in L^{p(\cdot)}(\Omega)$ and $g\in L^{p'(\cdot)}(\Omega)$ with a constant $r_{p(\cdot)}$ defined in \eqref{cr}.
\end{theorem}

\begin{definition}
 Given $p(\cdot)\in\mathcal{E}(\Omega)$, we define
\begin{align*}
&\||f\||_{p(\cdot)}=\sup_{m_{p'(\cdot)}(g)\le 1}\int_\Omega f(x)g(x)dx.
\end{align*}
\end{definition}

\begin{theorem}[equivalent norms, see \cite{KR}, Theorem 2.3]
Let $p(\cdot)\in \mathcal{E}(\Omega)$. Then the inequalities
\begin{align*}
&c_{p(\cdot)}^{-1}\|f\|_{p(\cdot)}\le \||f\||_{p(\cdot)}\le r_{p(\cdot)}\|f\|_{p(\cdot)}.
\end{align*}
 hold  for every $f\in L^{p(\cdot)}(\Omega)$, where the constants $c_{p(\cdot)}$ and  $r_{p(\cdot)}$   are  defined in \eqref{cr}.
\end{theorem}

We now define another norm and remind the reader that $\Omega$ is assumed to be an open bounded subset of $\mathbb{R}^N$.
\begin{definition}
Let $p\in\mathcal{E}(\Omega)$. Define
\begin{align*}
&]|f|[_{p(\cdot)}=\sup_{\|g\|_{p'(\cdot)}\le 1}\int_\Omega f(x)g(x)dx.
\end{align*}
\end{definition}

\begin{lemma}\label{dokcjocjrj}
Let $p\in\mathcal{E}(\Omega)$. Then for all $f\in L^{p(\cdot)}(\Omega)$,
\begin{align*}
&\||f\||_{p(\cdot)}=  ]|f|[_{p(\cdot)}.
\end{align*}
\end{lemma}
\begin{proof}
Set
\begin{align*}
&A=\{g;m_{p'(\cdot)}(g)\le 1\},\\
&B=\{g;\|g\|_{p'(\cdot)}\le 1\}.
\end{align*}
Clearly $m_{p'(\cdot)}(g)\le 1$ implies $\|g\|_{p'(\cdot)}\le 1$, which gives $A\subset B$. By the relation (2.11) in \cite{KR} we have
$m_{p'(\cdot)}(g)\le\|g\|_{p'(\cdot)}$ provided $\|g\|_{p'(\cdot)}\le 1$. Then $\|g\|_{p'(\cdot)}\le 1$ implies $m_{p'(\cdot)}(g)\le 1$ and so
$B\subset A$.

Thus, $A=B$ and
\begin{align*}
&\||f\||_{p(\cdot)}= \sup_{g\in A}\int_\Omega f(x)g(x)dx=\sup_{g\in B}\int_\Omega f(x)g(x)dx=]|f|[_{p(\cdot)}.
\end{align*}
\end{proof}
As a consequence we immediately obtain the following lemma.

\begin{lemma}\label{vjrfvi}
Let $p\in\mathcal{E}(\Omega)$. Then for all $f\in L^{p(\cdot)}(\Omega)$,
\begin{align*}
&c_{p(\cdot)}^{-1}\|f\|_{p(\cdot)}\le ]|f|[_{p(\cdot)}\le r_{p(\cdot)}\|f\|_{p(\cdot)},
\end{align*}
where the constants $c_{p(\cdot)}$ and  $r_{p(\cdot)}$   are  defined in \eqref{cr}.
\end{lemma}

\begin{theorem}[saturation of H\"older's inequality]\label{satur}
Let $p\in\mathcal{E}(\Omega)$. Then for each $\gamma<1$ and $f\in L^{p(\cdot)}(\Omega)$ there exists $g\in L^{p'(\cdot)}(\Omega)$ such that
\begin{align*}
&\|g\|_{p'(\cdot)}=1,\\
&\int_\Omega f(x)g(x)dx\ge \gamma\ c_{p(\cdot)}^{-1}\ \|f\|_{p(\cdot)}.
\end{align*}
\end{theorem}
\begin{proof}
Let $\gamma<1$ and $f\in L^{p(\cdot)}(\Omega)$ be given. If $f=0$ the proof is easy. Assume $f\neq 0$. Since
\begin{align*}
&]|f|[_{p(\cdot)}=\sup_{g\in B}\int_\Omega f(x)g(x)dx,
\end{align*}
where $B$ is given in the proof of Lemma \ref{dokcjocjrj}, there exists
 $g\in L^{p'(\cdot)}(\Omega)$ such that $\|g\|_{p'(\cdot)}\le 1$ and
\begin{align*}
&\int_\Omega f(x)g(x)dx\ge \gamma]|f|[_{p(\cdot)}.
\end{align*}
Clearly, $g\neq 0$.
Then by Lemma \ref{vjrfvi} we can write
\begin{align*}
&\int_\Omega f(x)g(x)dx\ge \gamma]|f|[_{p(\cdot)}\ge \gamma c_{p(\cdot)}^{-1}\|f\|_{p(\cdot)}\ge \gamma c_{p(\cdot)}^{-1}\|f\|_{p(\cdot)}\|g\|_{p'(\cdot)}.
\end{align*}
Now, take $g/{\|g\|_{p'(\cdot)}}$ instead of $g$ in the last inequality to obtain the assertion of our lemma.
\end{proof}

Given a function $u:\Omega\rightarrow\mathbb{R}$, its non-increasing rearrangement $u^*$   on $[0,\infty)$ is defined by
\begin{align*}
&u^*(t)=\inf\{\lambda>0;|\{x\in\Omega; \,\,\, |u(x)|>\lambda\}|\le t\}.
\end{align*}

\begin{lemma}\label{dcjdevnefovefjko}
Let $s:\Omega\rightarrow\mathbb{R}$ and $\alpha>1$. Then $(\alpha^{s(\cdot)})^*(t)=\alpha^{s^*(t)}$ for all $t> 0$.
\end{lemma}
\begin{proof}
Plainly
\begin{align*}
&(\alpha^{s(\cdot)})^*(t)=\inf\{\lambda>0;|\{x\in\Omega;\alpha^{s(x)}>\lambda\}|\le t\}\\
&=\inf\{\alpha^{\mu}>0;|\{x\in\Omega;\,\, \alpha^{s(x)}>\alpha^{\mu}\}|\le t\}\\
&=\inf\{\alpha^{\mu}>0;|\{x\in\Omega;\,\, {s(x)}>{\mu}\}|\le t\}\\
&=\alpha^{\inf\{{\mu}>0;|\{x\in\Omega;\,\, {s(x)}>{\mu}\}|\le t\}}=\alpha^{s^*(t)}.
\end{align*}
\end{proof}

In \cite{KR} (see Theorem 2.8) the following lemma is proved.
\begin{lemma}
Let   $p,\, q\in\mathcal{E}(\Omega)$. Then $L^{p(\cdot)}(\Omega)\hookrightarrow L^{q(\cdot)}(\Omega)$ if and only if $q(x)\le p(x)$ a.e. in $\Omega$.
\end{lemma}

\begin{definition}
A function $p:\Omega\rightarrow\mathbb{R}$ is said to satisfy a log-H\"older condition if there is $c>0$ such that
\begin{align}
&|p(x)- p(y)|\le -\frac{c}{\ln|x-y|},\ \ 0<|x-y|\le \frac12.\label{devjfjkovf}
\end{align}
\end{definition}

\begin{definition}
We say that $\Omega\in\mathcal{C}^{0,1}$ if there are a finite number of balls $\overline{B(x_k,r_k)},\, k=1,2,\dots,m$ and the same number of bi-Lipschitz mappings $T_k:[0,1]^{N-1}\times[-1,1]\rightarrow\overline{B(x_k,r_k)}$ such that for all $k\in\{1,2,\dots,m\}$,
\begin{enumerate}[\rm(i)]
\item $x_k\in\partial\Omega$,
\item $\bigcup_{k=1}^m B(x_k,r_k)\supset \partial\Omega$,
\item $T_k([0,1]^{N-1}\times[-1,0])= (\mathbb{R}^N\setminus \Omega)\cap \overline{B(x_k,r_k)}$,
\item $T_k([0,1]^{N-1}\times[0,1])=  \Omega\cap \overline{B(x_k,r_k)}$,
\item $T_k([0,1]^{N-1}\times\{0\})= \partial\Omega\cap \overline{B(x_k,r_k)}$.
\end{enumerate}
\end{definition}

Let $\Omega\in \mathcal{C}^{0,1}$ and $M\subset \overline{\Omega}$ be a compact set. Given $p(\cdot),\, q(\cdot):\overline{\Omega}\rightarrow\mathbb{R}$ we find conditions on $M$ and, moreover, determine how quickly can $q(\cdot)$ tend to $p^\#(\cdot)$ near $M$  while preserving the compactness of the embedding of $W^{1,p(\cdot)}(\Omega)$ in $L^{q(\cdot)}(\Omega)$.

\section{Almost-compact embedding between variable spaces}
Here we consider a set $\Omega\subset\mathbb{R}^N$
and functions $p(\cdot), \, q(\cdot)\in\mathcal{E}(\Omega)$; put $r(x)=p(x)/q(x)$. We will assume frequently that
\begin{align}
&\Omega\text{ is open and bounded,} \label{doicdocjdoij}\\
&q(x)\le p(x)\le p_+<\infty.\label{eovjoievjvj}
\end{align}

The assumption \eqref{doicdocjdoij} is assumed to be satisfied whenever variable exponent spaces are considered.

\begin{lemma}\label{sdcvdkovhiov}
Let  $q(\cdot)\in\mathcal{E}(\Omega)$, $q_{+}<\infty$. Then
\begin{align}
&\|u\|_{q(\cdot)}\le (m_{q(\cdot)}(u))^{1/{q_{+}}}\ \ \ \text{provided}\ \ \ m_{q(\cdot)}(u)\le 1,\label{devjkvnfjoebn}\\
&\|u\|_{q(\cdot)}\ge (m_{q(\cdot)}(u))^{1/{q_{-}}}\ \ \ \text{provided}\ \ \ m_{q(\cdot)}(u)\le 1,\label{jklfbjobno}\\
&\|u\|_{q(\cdot)}\le (m_{q(\cdot)}(u))^{1/{q_{-}}}\ \ \ \text{provided}\ \ \ m_{q(\cdot)}(u)\ge 1,\label{sdkcvbdovhbvhio}\\
&\|u\|_{q(\cdot)}\ge (m_{q(\cdot)}(u))^{1/{q_{+}}}\ \ \ \text{provided}\ \ \ m_{q(\cdot)}(u)\ge 1.\label{asdkcjdbvjhu}
\end{align}
\end{lemma}
\begin{proof}
Set $a=m_{q(\cdot)}(u)$ and assume $a\le 1$. Then
\begin{align*}
&\int_\Omega \Big(\frac{|u(x)|}{a^{1/q_{+}}}\Big)^{q(x)}dx\le\int_\Omega \Big(\frac{|u(x)|}{a^{1/q(x)}}\Big)^{q(x)}dx=\int_\Omega \frac{|u(x)|^{q(x)}}{a}dx=1
\end{align*}
which gives $\|u\|_{q(\cdot)}\le a^{1/q_{+}}$ and proves \eqref{devjkvnfjoebn}. The assertions \eqref{jklfbjobno}, \eqref{sdkcvbdovhbvhio} and \eqref{asdkcjdbvjhu} can be proved analogously.
\end{proof}

\begin{lemma}\label{sddovhvhbn}
Assume that $\Omega,\,  p(\cdot),\, q(\cdot)$ satisfy \eqref{doicdocjdoij} and \eqref{eovjoievjvj}. Let  $\|u\|_{p(\cdot)}\le 1$. Then
\begin{align*}
&\|u\|_{p(\cdot)}^{q_{+}}\le\|\ |u(\cdot)|^{q(\cdot)}\|_{r(\cdot)}\le \|u\|_{p(\cdot)}^{q_{-}}.
\end{align*}
\end{lemma}
\begin{proof}
Assume first $0<a:=\|u\|_{p(\cdot)}<1$. Hence
\begin{align}
&1=\int_\Omega \Big(\frac{|u(x)|}{a}\Big)^{p(x)}dx>\int_\Omega |u(x)|^{p(x)}dx.\label{sjocdjcvhdwvh}
\end{align}
Set $b=\|\ |u(\cdot)|^{q(\cdot)}\|_{r(\cdot)}$. Then
\begin{align*}
&1=\int_\Omega \Big(\frac{|u(x)|^{q(x)}}{b}\Big)^{r(x)}dx=\int_\Omega \Big(\frac{|u(x)|}{b^{1/{q(x)}}}\Big)^{p(x)}dx.
\end{align*}
If $b>1$, then
\begin{align*}
&1=\int_\Omega \Big(\frac{|u(x)|}{b^{1/{q(x)}}}\Big)^{p(x)}dx<\int_\Omega |u(x)|^{p(x)}dx \overset{\eqref{sjocdjcvhdwvh}}{<}1,
\end{align*}
which is a contradiction. So $b\le 1$. Consequently,
\begin{align*}
&\int_\Omega\Big(\frac{|u(x)|}{b^{1/{q_{+}}}}\Big)^{p(x)}dx\le 1=\int_\Omega \Big(\frac{|u(x)|}{b^{1/{q(x)}}}\Big)^{p(x)}dx\le\int_\Omega \Big(\frac{|u(x)|}{b^{1/{q_{-}}}}\Big)^{p(x)}dx,
\end{align*}
which gives $b^{1/q_{+}}\ge \|u\|_{p(\cdot)}\ge b^{1/q_{-}}$ and finally
\begin{align*}
&\|u\|_{p(\cdot)}^{q_{+}}\le \|\ |u(\cdot)|^{q(\cdot)}\|_{r(\cdot)}\le\|u\|_{p(\cdot)}^{q_{-}}.
\end{align*}

Assume now $\|u\|_{p(\cdot)}=1$. Choose $\varepsilon>0$. Then $\|\frac{u}{1+\varepsilon}\|_{p(\cdot)}<1$ and so,
\begin{align*}
&\Big\|\ \frac{|u(\cdot)|}{1+\varepsilon}\Big\|_{{p(\cdot)}}^{q_{+}} \le \Big\|\ \Big(\frac{|u(\cdot)|}{1+\varepsilon}\Big)^{q(\cdot)}\Big\|_{r(\cdot)}\le \Big\|\ \frac{|u(\cdot)|}{1+\varepsilon}\Big\|_{{p(\cdot)}}^{q_{-}}.
\end{align*}
Since $\frac{1}{(1+\varepsilon)^{q_{-}}}\ge \frac{1}{(1+\varepsilon)^{q(x)}}\ge \frac{1}{(1+\varepsilon)^{q_{+}}}$ we have
\begin{align*}
&\frac{1}{(1+\varepsilon)^{q_{+}}}\Big\|\ |u(\cdot)|^{q(\cdot)}\Big\|_{r(\cdot)}\le\Big\|\ \Big(\frac{|u(\cdot)|}{1+\varepsilon}\Big)^{q(\cdot)}\Big\|_{r(\cdot)}\le
\Big\|\ \frac{|u(\cdot)|}{1+\varepsilon}\Big\|_{p(\cdot)}^{q_{-}}
\le \frac{1}{(1+\varepsilon)^{q_{-}}}\|u\|_{p(\cdot)}^{q_{-}}
\end{align*}
and
\begin{align*}
&\frac{1}{(1+\varepsilon)^{q_{-}}}\Big\|\ |u(\cdot)|^{q(\cdot)}\Big\|_{r(\cdot)}\ge\Big\|\ \Big(\frac{|u(\cdot)|}{1+\varepsilon}\Big)^{q(\cdot)}\Big\|_{r(\cdot)}\ge
\Big\|\ \frac{|u(\cdot)|}{1+\varepsilon}\Big\|_{p(\cdot)}^{q_{+}}
\ge \frac{1}{(1+\varepsilon)^{q_{+}}}\|u\|_{p(\cdot)}^{q_{+}},
\end{align*}
which proves
\begin{align*}
&(1+\varepsilon)^{q_{-}-q_{+}}\|u\|_{p(\cdot)}^{q_{+}}\le \|\ |u(\cdot)|^{q(\cdot)}\|_{r(\cdot)}\le(1+\varepsilon)^{q_{+}-q_{-}}\|u\|_{p(\cdot)}^{q_{-}}.
\end{align*}
Letting $\varepsilon\rightarrow 0^+$ we conclude that
\begin{align*}
&\|u\|_{p(\cdot)}^{q_{+}}\le \|\ |u(\cdot)|^{q(\cdot)}\|_{r(\cdot)}\le\|u\|_{p(\cdot)}^{q_{-}}.
\end{align*}
\end{proof}

\begin{lemma}\label{skovkvrbn}Assume that $\Omega, \, p(\cdot), \, q(\cdot)$ satisfy \eqref{doicdocjdoij} and \eqref{eovjoievjvj}. Suppose also that for any sequence $\{E_n\}_{n\in \N}$ of measurable subsets of $\Omega$ such that $|E_n|\rightarrow 0$, we have
\begin{align*}
&\|\chi_{E}\|_{r'(\cdot)}\rightarrow 0.
\end{align*}
Then $L^{p(\cdot)}(\Omega)\overset{\ast}{\hookrightarrow} L^{q(\cdot)}(\Omega)$.
\end{lemma}
\begin{proof}
Let $E_n\subset \Omega$, $|E_n|\rightarrow 0$. Then by Lemma \ref{sdcvdkovhiov} we obtain
\begin{align*}
&\lim_{n\rightarrow\infty}\sup\{\|u\chi_{E_n}\|_{q(\cdot)};\, \|u\|_{p(\cdot)}\le 1\}\\
&\le \lim_{n\rightarrow\infty}\sup\{\max\{(m_{q(\cdot)}(u\chi_{E_n}))^{1/{q_{+}}},(m_{q(\cdot)}(u\chi_{E_n}))^{1/{q_{-}}}\};\, \|u\|_{p(\cdot)}\le 1\}.
\end{align*}
If $\|u\|_{p(\cdot)}\le 1$, we obtain by the H\"older inequality and Lemma \ref{sddovhvhbn}
\begin{align*}
&m_{q(\cdot)}(u\chi_{E_n})=\int_\Omega |u(x)\chi_{E_n}(x)|^{q(x)}dx\le c\|\chi_{E_n}\|_{r'(\cdot)}\|\ |u(\cdot)|^{q(\cdot)}\|_{r(\cdot)}\\
&\le c\|\chi_{E_n}\|_{r'(\cdot)}\|u\|_{p(\cdot)}^{q_{-}}.
\end{align*}
This gives
\begin{align*}
&\lim_{n\rightarrow\infty}\sup\{\|u\chi_{E_n}\|_{q(\cdot)};\, \|u\|_{p(\cdot)}\le 1\}\\
&\le c\lim_{n\rightarrow\infty}\sup\{\max\{\|\chi_{E_n}\|_{r'(\cdot)}^{1/{q_{+}}}\|u\|_{p(\cdot)}^{q_{-}/{q_{+}}},\, \|\chi_{E_n}\|_{r'(\cdot)}^{1/{q_{-}}}\|u\|_{p(\cdot)}\};\, \|u\|_{p(\cdot)}\le 1\}\\
&\le c\lim_{n\rightarrow\infty}\max\{ \|\chi_{E_n}\|_{r'(\cdot)}^{1/{q_{+}}},\, \| \chi_{E_n}\|_{r'(\cdot)}^{1/{q_{-}}}\}=0.
\end{align*}

\end{proof}

Put for each $x\in\Omega$
\begin{align*}
  s(x)=
  \begin{cases}
  \frac{1}{p(x)-q(x)}, &\ p(x)>q(x),\\
  \infty, &\ p(x)=q(x).
  \end{cases}
\end{align*}

\begin{theorem} \label{thm3.4}
Assume that $\Omega, \, p(\cdot),\, q(\cdot)$ satisfy \eqref{doicdocjdoij} and \eqref{eovjoievjvj}.
 Assume
\begin{align}
&\int_{0}^{|\Omega|}a^{s^\ast(t)}dt<\infty\label{sdjcdjvfjo}
\end{align}
for all $a>1$ (we set $a^\infty=\infty$). Then $L^{p(\cdot)}(\Omega)\overset{\ast}{\hookrightarrow} L^{q(\cdot)}(\Omega)$.
\end{theorem}
\begin{proof}
Let $E_n\subset \Omega$, $|E_n|\rightarrow 0$. Assume that there is $\alpha>0$ such that
\begin{align*}
&\|\chi_{E_n}\|_{r'(\cdot)}\ge \alpha
\end{align*}
for all $n$. Without loss of generality we can suppose $\alpha<1$. Then
\begin{align*}
&\alpha\le \inf\Big\{\lambda>0;\int_\Omega \Big|\frac{\chi_{E_n}(x)}{\lambda}\Big|^{r'(x)}dx\le 1\Big\}= \inf\Big\{\lambda>0;\int_\Omega \Big|\frac{\chi_{E_n}(x)}{\lambda}\Big|^{p(x)s(x)}dx\le 1\Big\}.
\end{align*}
Choose $0<\beta<\alpha$. Then we obtain by Lemma \ref{dcjdevnefovefjko},
\begin{align*}
&1< \int_\Omega \Big|\frac{\chi_{E_n}(x)}{\beta}\Big|^{p(x)s(x)}dx\le\int_{E_n} \Big(\frac{1}{\beta^{p_{+}}}\Big)^{s(x)}dx\le\int_0^{|E_n|}\Big(\frac{1}{\beta^{p_{+}}}\Big)^{s^*(t)}dt.
\end{align*}
Since by the assumption,
\begin{align*}
&\int_0^{|E_n|}\Big(\frac{1}{\beta^{p_{+}}}\Big)^{s^*(t)}dt \rightarrow 0 \ \ \text{for}\ \ n\rightarrow\infty,
\end{align*}
we have a contradiction. Then for each subsequence $E_{n_k}$, $|E_{n_k}|\rightarrow 0$, we have $\liminf_{k\rightarrow\infty}\|\chi_{E_{n_k}}\|_{r'(\cdot)}=0$. So, $\|\chi_{E_n}\|_{r'(\cdot)}\rightarrow 0$. By Lemma \ref{skovkvrbn} we have $L^{p(\cdot)}(\Omega)\overset{\ast}{\hookrightarrow} L^{q(\cdot)}(\Omega)$.
\end{proof}

Note that if \eqref{sdjcdjvfjo} is satisfied it is easy to see that $p\not=q$ almost everywhere.
We remark that the condition \eqref{sdjcdjvfjo} was first introduced in Corollary 2.7 of \cite{ELN}.

\begin{lemma}\label{devdevhi}
Assume that $\Omega,\,  p(\cdot),\, q(\cdot)$ satisfy \eqref{doicdocjdoij} and \eqref{eovjoievjvj}. Assume that there exist an $\alpha>0$ and a sequence $E_n\subset \Omega$ with $|E_n|\rightarrow 0$ such that for all $n$
\begin{align*}
&\|\chi_{E_n}\|_{r'(\cdot)}\ge \alpha.
\end{align*}
Then $L^{p(\cdot)}(\Omega)\overset{\ \ast}{\hookrightarrow\!\!\!\!\!\!\!\text{\scriptsize{/}}}\ L^{q(\cdot)}(\Omega)$.
\end{lemma}
\begin{proof}
Let $E_n\subset \Omega$, $\|\chi_{E_n}\|_{r'(\cdot)}\ge \alpha$. Fix $n$. Without loss of generality we can assume $\al \le 1$. By Theorem \ref{satur} there exists $g_n\in L^{r(\cdot)}(\Omega)$ such that\begin{align*}
&\|g_n\|_{r(\cdot)}=1,\\
&\int_\Omega g_n(x)\chi_{E_n}(x)dx\ge \frac{1}{2}\ c_{r'(\cdot)}^{-1}\ \|\chi_{E_n}\|_{r'(\cdot)}.
\end{align*}
Set $u_n(x)=g_n(x)^{1/{q(x)}}$.
Clearly, since $\|g_n\|_{r(\cdot)}=1$ we obtain by (2.11) in \cite{KR},
\begin{align*}
&\int_\Omega |u_n(x)|^{p(x)}dx=\int_\Omega g_n(x)^{r(x)} dx=m_{r(\cdot)}(g_n)\le \|g_n\|_{r(\cdot)}=1.
\end{align*}
Thus
\begin{align*}
&\|u_n\|_{p(\cdot)}\le 1.
\end{align*}

Further
\begin{align*}
&m_{q(\cdot)}(u_n\chi_{E_n})=\int_\Omega |u_n(x)\chi_{E_n}(x)|^{q(x)}dx=\int_\Omega |g_n(x)|\chi_{E_n}(x)dx\\
&\ge \int_\Omega g_n(x)\chi_{E_n}(x)dx\ge\frac{1}{2}\ c_{r'(\cdot)}^{-1}\ \|\chi_{E_n}\|_{r'(\cdot)}\ge \frac{1}{2}\ c_{r'(\cdot)}^{-1}\ \alpha=:\beta.
\end{align*}
By Lemma \ref{sdcvdkovhiov} we obtain
\begin{align*}
\sup\{\|u\chi_{E_n}\|_{q(\cdot)};\, \|u\|_{p(\cdot)}\le 1\}&\ge \min\{(m_{q(\cdot)}(u_n\chi_{E_n}))^{1/{q_{+}}},\, (m_{q(\cdot)}(u_n\chi_{E_n}))^{1/{q_{-}}}\}\\
&\ge \min\{\beta^{1/q_{+}},\beta^{1/q_{-}}\}=\beta^{1/q_{-}},
\end{align*}
which completes the proof.
\end{proof}

\begin{lemma}\label{bkobhnrhjknmpk}
Let $E\subset\Omega$, $g:\Omega\rightarrow[0,\infty)$ and assume that
\begin{align*}
&\inf\{g(x);x\in E\}\ge\sup\{g(x);x\in \Omega\setminus E\}.
\end{align*}
Then
\begin{align*}
&(g\chi_E)^*(t)=g^*(t)\chi_{(0,|E|)}(t).
\end{align*}
\end{lemma}
\begin{proof}
As the proof is trivial it is omitted.
\end{proof}

\begin{theorem}
Assume that $\Omega,\,  p(\cdot),\, q(\cdot)$ satisfy \eqref{doicdocjdoij} and \eqref{eovjoievjvj}. Suppose that there is $a>1$ such that
\begin{align}
&\int_{0}^{|\Omega|}a^{s^\ast(t)}dt=\infty.\label{rjvtortbb}
\end{align}

 Then $L^{p(\cdot)}(\Omega)\overset{\ \ast}{\hookrightarrow\!\!\!\!\!\!\!\text{\scriptsize{/}}}\ L^{q(\cdot)}(\Omega)$.
\end{theorem}
\begin{proof}
Define $E_n=\{x\in\Omega;s(x)\ge n\}$. Assume that there exists $n_0$ such that $|E_{n_0}|=0$. Then
\begin{align*}
&s(x)=\frac{1}{p(x)-q(x)}<n_0
\end{align*}
almost everywhere and so we have for any $a>1$,
\begin{align*}
&\int_0^{|\Omega|}a^{s^*(t)}dt\le\int_0^{|\Omega|}a^{n_0}dt=a^{n_0}|\Omega|<\infty,
\end{align*}
which contradicts the assumption. So, $|E_n|>0$ for all $n$. Fix $n$ and assume
\begin{align}
&\max\{\|\chi_{E_n}\|_{r'(\cdot)}^{p_{+}},\, \|\chi_{E_n}\|_{r'(\cdot)}^{p_{-}}\}\le \frac{1}{a}.\label{devfjklebvnefn}
\end{align}
Then
\begin{align}
&1\ge\int_\Omega \left(\frac{\chi_{E_n}(x)}{\|\chi_E\|_{r'(\cdot)}}\right)^{r'(x)}dx=\int_\Omega \left(\frac{\chi_{E_n}(x)}{\|\chi_{E_n}\|_{r'(\cdot)}^{p(x)}}\right)^{s(x)}dx\label{cvvjioobvhiobip}\\
&\ge \int_{E_n}\left(\frac{1}{\max\{\|\chi_{E_n}\|_{r'(\cdot)}^{p_{+}},\|\chi_{E_n}\|_{r'(\cdot)}^{p_{-}}\}}\right)^{s(x)}dx\ge \int_{E_n} a^{s(x)}dx.\nonumber
\end{align}
Now, by the definition of $E_n$ we have that $s(x)\ge n$ on $E_n$ and $s(x)<n$ on $\Omega\setminus E_n$.  This gives us $a^{s(x)}\ge a^n$ on $E_n$ and $a^{s(x)}< a^n$ on $\Omega\setminus E_n$.  Then we have by Lemma \ref{bkobhnrhjknmpk}
\begin{align*}
&(a^{s(\cdot)}\chi_{E_n}(\cdot))^*(t)=(a^{s(\cdot)})^*(t)\chi_{{(0,|E_n|)}}(t)=a^{s^*(t)}\chi_{{(0,|E_n|)}}(t),
\end{align*}
which gives with \eqref{cvvjioobvhiobip}
\begin{align*}
&1\ge \int_{E_n} a^{s(x)}dx=\int_\Omega a^{s(x)}\chi_{E_n}(x)dx=\int_0^{|\Omega|}(a^{s(\cdot)}\chi_{E_n}(\cdot))^*(t)dt\\
&=\int_0^{|\Omega|}a^{s^*(t)}\chi_{{(0,|E_n|)}}(t)dt=\int_0^{|E_n|}a^{s^*(t)}(t)dt=\infty,
\end{align*}
which is a contradiction.
Hence our assumption \eqref{devfjklebvnefn} is false and we have
\begin{align*}
&\max\{\|\chi_{E_n}\|_{r'(\cdot)}^{p_{+}},\|\chi_{E_n}\|_{r'(\cdot)}^{p_{-}}\}> \frac{1}{a},
\end{align*}
which yields
\begin{align*}
&\|\chi_{E_n}\|_{r'(\cdot)}> \min\{a^{-1/{p_{+}}}  ,a^{-1/{p_{-}}}\}=:b>0.
\end{align*}
Thus, we have $\|\chi_{E_n}\|_{r'(\cdot)}\ge b>0$ for any $n$ and Lemma \ref{devdevhi} gives us $L^{p(\cdot)}(\Omega)\overset{\ \ast}{\hookrightarrow\!\!\!\!\!\!\!\text{\scriptsize{/}}}\ L^{q(\cdot)}(\Omega)$.
\end{proof}

Now suppose that $K\subset\Omega$ is compact with $|K|=0$. Denote $d_K(x)=\dist(x,K)$, set
\begin{align}
&K(t)=\{x\in\Omega;\d_K(x)<t\},\label{jkveebhjiopbjtiopb}
\end{align}
and put
\begin{align}
&\varphi(t)=|K(t)|,\,\, t\in [0,\diam(\Omega)].\label{efovhevhiobh}
\end{align}
Let $\omega:(0,\diam(\Omega)]\rightarrow\mathbb{R}$ be a decreasing continuous non-negative function,
$\omega_0:=\omega(\diam(\Omega))$. We set $\omega(0)=\lim_{t\rightarrow 0_+}\omega(t)$.
Let $\omega^{-1}$ denote the inverse function to $\omega$.

\begin{lemma}\label{fkvjkkvjovji}
Assume that $\Omega,\,  p(\cdot),\, q(\cdot)$ satisfy \eqref{doicdocjdoij} and \eqref{eovjoievjvj}. Assume that there is $c>0$ such that
\begin{align*}
&s(x)=\frac{1}{p(x)-q(x)}\le c\ \omega(\d_K(x)), \ x\in \Omega,\\
&\int_{\omega_0}^{\omega(0)} \varphi(\omega^{-1}(y))a^y dy<\infty \ \text{for all}\ a>1.
\end{align*}
Then $L^{p(\cdot)}(\Omega)\overset{\ \ast}{\hookrightarrow} L^{q(\cdot)}(\Omega)$.
\end{lemma}
\begin{proof}

 Let $a>1$. Then
\begin{align*}
&\int_0^{|\Omega|}a^{s^*(t)}dt=\int_\Omega a^{s(x)}dx\le\int_\Omega a^{c\omega(\d_K(x))}dx=
\int_0^{\infty} |\{x;a^{c\omega(\d_K(x))}>\lambda\}|d\lambda\\
&=\int_0^{a^{c\omega_0}} |\{x;a^{c\omega(\d_K(x))}>\lambda\}|d\lambda+\int_{a^{c\omega_0}} ^{\infty} |\{x;a^{c\omega(\d_K(x))}>\lambda\}|d\lambda\\
&=\int_0^{a^{c\omega_0}} |\{x;a^{c\omega(\d_K(x))}>\lambda\}|d\lambda+\int_{a^{c\omega_0}} ^{a^{c\omega(0)}} |\{x;a^{c\omega(\d_K(x))}>\lambda\}|d\lambda\\
&\qquad\qquad                                 +\int_{a^{c\omega(0)}} ^{\infty} |\{x;a^{c\omega(\d_K(x))}>\lambda\}|d\lambda.
\end{align*}
Since
\begin{align*}
&\int_0^{a^{c\omega_0}} |\{x;a^{c\omega(\d_K(x))}>\lambda\}|d\lambda={a^{c\omega_0}} |\Omega|
\end{align*}
and
\begin{align*}
&|\{x;a^{c\omega(\d_K(x))}>\lambda\}|=0\text{ for }\lambda>a^{c\omega(0)},
\end{align*}
we obtain
\begin{align*}
&\int_0^{|\Omega|}a^{s^*(t)}dt={a^{c\omega_0}} |\Omega|+\int_{a^{c\omega_0}} ^{a^{c\omega(0)}} |\{x;a^{c\omega(\d_K(x))}>\lambda\}|d\lambda\\
&={a^{c\omega_0}} |\Omega|+c\ln a\int_{\omega_0} ^{\omega(0)} |\{x;a^{c\omega(\d_K(x))}>a^{cy}\}|a^{cy}dy\\
&={a^{c\omega_0}} |\Omega|+c\ln a\int_{\omega_0} ^{\omega(0)} |\{x;{\omega(\d_K(x))}>{y}\}|a^{cy}dy\\
&={a^{c\omega_0}} |\Omega|+c\ln a\int_{\omega_0} ^{\omega(0)} |\{x;\d_K(x)<\omega^{-1}(y)\}|a^{cy}dy\\
&={a^{c\omega_0}} |\Omega|+c\ln a\int_{\omega_0} ^{\omega(0)} \varphi(\omega^{-1}(y))a^{cy}dy<\infty.
\end{align*}
By Theorem~\ref{thm3.4} we have $L^{p(\cdot)}(\Omega)\overset{\ast}{\hookrightarrow} L^{q(\cdot)}(\Omega)$.
\end{proof}

\section{Examples of Cantor sets}

Let $\{a_k\}_{k\in \N}$ be a given sequence of positive real numbers with
\begin{align}
&\sum_{k=1}^\infty a_k=1.\label{vjklfjkovnkobo}
\end{align}
Construct a generalized Cantor set by the following process. Set $K_0=[0,1]$. Omit in the first step from $K_0$ a centered interval of length $a_1$ to obtain a set $K_1$. We write
\begin{align*}
&K_1=K_0\setminus \Big(\frac{1-a_1}{2},\frac{1+a_1}{2}\Big)=\Big[0,\frac{1-a_1}{2}\Big]\cup\Big[\frac{1+a_1}{2},1\Big]:=J_0\cup J_1.
\end{align*}
In the second step we omit from $J_0$ and $J_1$ centered intervals of length $a_2/2$ to obtain $K_2$. Then
\begin{align*}
&K_2=K_1\setminus \Big(\Big(\frac{1-a_1-a_2}{2^2},\frac{1-a_1+a_2}{2^2}\Big)\cup\Big(\frac{3+a_1-a_2}{2^2},\frac{3+a_1+a_2}{2^2}\Big)\Big)\\
&\!=\!\Big[0,\frac{1-a_1-a_2}{2^2}\Big]\!\cup\!\Big[\frac{1-a_1+a_2}{2^2},\frac{1-a_1}{2}\Big]\\
&\qquad\!\cup\! \Big[\frac{1+a_1}{2},\frac{3+a_1-a_2}{2^2}\Big]\!\cup\!\Big[\frac{3+a_1+a_2}{2^2},1\Big]:=J_{00}\cup J_{01}\cup J_{10}\cup J_{11}.
\end{align*}
We follow this process step by step to obtain sets $K_n$. Then $K_n$ consists of $2^n$ intervals $J_\alpha$, $\alpha\in\{0,1\}^n$. Clearly,
\begin{align}
&|J_\alpha|=2^{-n}\Big(1-\sum_{k=1}^n a_k\Big). \label{djvdjklvfrbv}
\end{align}
Set
\begin{align*}
&K=\bigcap_{n=1}^\infty K_n.
\end{align*}
Clearly, $K$ is a compact set and for each $n$
\begin{align*}
&|K|\le |K_n|,
\end{align*}
which gives with \eqref{djvdjklvfrbv}
\begin{align*}
&|K|\le \sum_{\alpha\in\{0,1\}^n}|J_\alpha|=2^n 2^{-n}\Big(1-\sum_{k=1}^n a_k\Big)=1-\sum_{k=1}^n a_k.
\end{align*}
Using \eqref{vjklfjkovnkobo} we have
\begin{align*}
&|K|=0.
\end{align*}

Now, we will be interested in the behavior of the function $|K(t)|$ defined by \eqref{jkveebhjiopbjtiopb}.
\begin{lemma}\label{dwjkocvdovrov}
The function $|K(\cdot)|$ is non-decreasing and $\lim_{t\rightarrow 0_+}|K(t)|=0$.
\end{lemma}
\begin{proof}
The monotonicity of $|K(\cdot)|$ is clear. Moreover, $K(t)\searrow K$ and $K(1)<\infty$ since $K$ is compact. It is easily seen that $\lim_{t\rightarrow 0_+}|K(t)|=0$.
\end{proof}
\begin{lemma}\label{decvefjvjjo}
For each $n\in\mathbb{N}$ let $r_n, \varepsilon_n$ be given by
\begin{align*}
&r_n=1-\sum_{k=1}^n a_k,\\
&\varepsilon_n=2^{-n}\Big(1-\sum_{k=1}^n a_k\Big)=2^{-n}r_n.
\end{align*}
Then
\begin{align*}
&r_n\le |K(\varepsilon_n)|\le 4r_n.
\end{align*}
\end{lemma}
\begin{proof}
Clearly, $\varepsilon_n=|J_\alpha|$ for $\alpha\in\{0,1\}^n$ and so,
\begin{align*}
&K(\varepsilon_n)\supset \bigcup_{\alpha\in\{0,1\}^n} J_{\alpha},
\end{align*}
which gives
\begin{align*}
&|K(\varepsilon_n)|\ge \sum_{\alpha\in\{0,1\}^n} |J_{\alpha}|=2^n\varepsilon_n=r_n.
\end{align*}
For the right-hand inequality, denote by $M$ the set of all endpoints of intervals $J_{\alpha}$, $\alpha\in\{0,1\}^n$. The number of these points is $2(1+2+2^2+\dots+2^{n-1})=2(2^n-1)$ and
\begin{align*}
&K(\varepsilon_n)\subset\bigcup_{x\in M}(x-\varepsilon_n,x+\varepsilon_n).
\end{align*}
Then
\begin{align*}
&|K(\varepsilon_n)|\le \sum_{x\in M}2\varepsilon_n=2(2^n-1)2\varepsilon_n\le 4\cdot 2^n\varepsilon_n=4r_n.
\end{align*}
\end{proof}

One important case is obtained by choosing $a_k=\frac{a^{k-1}}{(a+1)^{k}}$ where $a>0$. When $a=2$ we obtain the classical Cantor set.
\begin{lemma}\label{djkojkovnfovn}
Let $a_k=\frac{a^{k-1}}{(a+1)^{k}}$ and set
\begin{align*}
&s=\frac{\ln(\frac{a}{a+1})}{\ln(\frac{a}{2(a+1)})}.
\end{align*}
Then there are positive constants $c_1,c_2$ such that
\begin{align*}
&c_1 t^s\le |K(t)|\le c_2 t^s,\ t\in [0,1].
\end{align*}
\end{lemma}
\begin{proof}
Let $q=\frac{a}{a+1}$. Then $s=\frac{\ln q}{\ln(q/2)}$. Let $r_n$ and $\varepsilon_n$ be as defined in Lemma~\ref{decvefjvjjo}. Clearly,
\begin{align*}
&r_n=1-\sum_{k=1}^{n}\frac{a^{k-1}}{(a+1)^{k}}=\Big(\frac{a}{a+1}\Big)^n=q^n,\\
&\varepsilon_n=2^{-n}\Big(\frac{a}{a+1}\Big)^n=2^{-n}q^n
\end{align*}
and
\begin{align*}
&r_{n+1}=qr_n,\ \ \  \varepsilon_{n+1}=\frac{q}{2}\varepsilon_n.
\end{align*}
It is easy to see that $0<s<1$.

Fix $t\in[\varepsilon_{n+1},\varepsilon_{n}]$. By Lemmas \ref{dwjkocvdovrov} and \ref{decvefjvjjo} we know that
\begin{align}
&qr_n=r_{n+1}\le |K(\varepsilon_{n+1})|\le |K(t)|\le |K(\varepsilon_{n})|\le 4r_n.\label{cvhbrojverjobv}
\end{align}
Since $(q/2)\varepsilon_{n}\le t\le \varepsilon_{n}$, we have
\begin{align*}
&\ln (q/2)+n\ln (q/2)=\ln (q/2)+\ln \varepsilon_{n}\le \ln t\le \ln\varepsilon_{n}=n\ln (q/2),
\end{align*}
which gives
\begin{align*}
&\frac{\ln t}{\ln (q/2)}\ge n\ge \frac{\ln t-\ln (q/2)}{\ln (q/2)}.
\end{align*}
This implies that
\begin{align*}
&t^{\frac{\ln q}{\ln (q/2)}}=q^{\frac{\ln t}{\ln (q/2)}}\le r_n\le q^{\frac{\ln t-\ln (q/2)}{\ln (q/2)}}=\frac{1}{q}\ t^{\frac{\ln q}{\ln (q/2)}}.
\end{align*}
By \eqref{cvhbrojverjobv} we obtain
\begin{align*}
&q t^s\le |K(t)|\le \frac{4}{q}\ t^s.
\end{align*}
\end{proof}

We recall the definition of the Riemann function
\begin{align*}
&\zeta(s)=\sum_{k=1}^\infty k^{-s}, \ s>1.
\end{align*}
\begin{lemma}\label{vjvjbvjvvjvjj}
Let $s>1$ and  $a_k=\frac{k^{-s}}{\zeta(s)}$.
Then there are positive constants $c_1,c_2$ such that
\begin{align*}
&\frac{c_1}{(\ln(\e/t))^{s-1}} \le |K(t)|\le \frac{c_2}{(\ln(\e/t))^{s-1}},\ t\in [0,1].
\end{align*}
\end{lemma}
\begin{proof} Let $r_n$ and $\varepsilon_n$ be as  defined in Lemma~\ref{decvefjvjjo}.
 Clearly,
\begin{align*}
&r_n=1-\frac{1}{\zeta(s)}\sum_{k=1}^{n}k^{-s}=\frac{1}{\zeta(s)}\sum_{k=n+1}^{\infty}k^{-s},\ \varepsilon_n=2^{-n}r_n.
\end{align*}
It is easy to see that
\begin{align*}
&\frac{1}{\zeta(s)(s-1)(n+1)^{s-1}}=\frac{1}{\zeta(s)}\int_{n+1}^\infty x^{-s}ds\le r_n\\
&\le \frac{1}{\zeta(s)}\int_{n}^\infty x^{-s}ds=\frac{1}{\zeta(s)(s-1)n^{s-1}}.
\end{align*}
This gives for $ n\ge 2$,
\begin{align*}
&2^{1-s}\le \frac{r_{n+1}}{r_n}\le 1,\ \ \ \ 2^{-s}\le \frac{\varepsilon_{n+1}}{\varepsilon_n}\le \frac{1}{2}.
\end{align*}

Fix $t\in[\varepsilon_{n+1},\varepsilon_{n}]$. Then
\begin{align}
&2^{1-s} r_n\le r_{n+1}\le |K(\varepsilon_{n+1})|\le |K(t)|\le |K(\varepsilon_{n})|\le 4r_n.\label{efvfobgkibnnbko}
\end{align}
We know
\begin{align*}
&2^{-s}\varepsilon_{n}\le \varepsilon_{n+1}\le t\le \varepsilon_n,
\end{align*}
which gives
\begin{align*}
&\frac{2^{-n}(n+1)^{1-s}}{\zeta(s)(s-1)}\le \varepsilon_n\le \frac{2^{-n}n^{1-s}}{\zeta(s)(s-1)}
\end{align*}
and consequently
\begin{align*}
&\frac{1}{2^s2^{s-1}}\frac{2^{-n}n^{1-s}}{\zeta(s)(s-1)}\le t\le \frac{2^{-n}n^{1-s}}{\zeta(s)(s-1)}.
\end{align*}
So, there are constants $b_1, \,b_2$ such that
\begin{align*}
&b_1 2^{-n} n^{1-s}\le t\le b_2 2^{-n} n^{1-s}.
\end{align*}
Hence
\begin{align*}
&\ln b_1-n\ln 2-(s-1) n\le\ln b_1-n\ln 2-(s-1)\ln n\\
&\qquad\qquad\le \ln t\\
&\le \ln b_2-n\ln 2-(s-1)\ln n\le \ln b_2-n\ln 2
\end{align*}
and thus
\begin{align*}
&\frac{\ln(b_1/t)}{\ln 2+s-1}\le n\le \frac{\ln(b_2/t)}{\ln 2}.
\end{align*}
Then
\begin{align*}
&|K(t)|\le 4r_n\le\frac{4}{\zeta(s)(s-1)n^{s-1}}\le\frac{4(\ln 2+s-1)^{s-1}}{\zeta(s)(s-1)(\ln(b_1/t))^{s-1}}\le\frac{c_2}{(\ln(\e/t))^{s-1}}.
\end{align*}
Finally,
\begin{align*}
&|K(t)|\ge r_{n+1}\ge\frac{4}{\zeta(s)(s-1)(n+2)^{s-1}}\ge\frac{1}{\zeta(s)(s-1)}\Big(\frac{n}{n+2}\Big)^{s-1}\frac{1}{n^{s-1}}\\
&\ge \frac{1}{3^{s-1}\zeta(s)(s-1)n^{s-1}}\ge\frac{(\ln 2)^{s-1}}{3^{s-1}\zeta(s)(s-1)(\ln(b_2/t))^{s-1}}\ge\frac{c_1}{(\ln(\e/t))^{s-1}}.
\end{align*}
\end{proof}

\begin{lemma}\label{dsjkvbdjovdjvj}
Define a function $\eta(s)=\sum_{k=1}^\infty \frac{1}{(k+1)\ln^s(k+1)}$, $s>1$. Choose $a_k=\frac{1}{\eta(s)}\frac{1}{(k+1)\ln^s(k+1)}$.
Then there are positive constants $c_1,\, c_2$ and $b$ such that
\begin{align*}
&c_1 (\ln\ln(b/t))^{1-s}\le |K(t)|\le c_2  (\ln\ln(b/t))^{1-s},\ t\in[0,1].
\end{align*}
\end{lemma}
\begin{proof} The proof is analogous to the previous one. Clearly,
\begin{align*}
&r_n=\frac{1}{\eta(s)}\sum_{k=n+1}^\infty \frac{1}{(k+1)\ln^s(k+1)},\ \ \varepsilon_n=2^{-n}r_n.
\end{align*}
By the integral criterion we  have the estimates 
\begin{align*}
&\frac{1}{\eta(s)(s-1)\ln^{s-1}(n+2)}\le r_n\le\frac{1}{\eta(s)(s-1)\ln^{s-1}(n+1)}, \quad n\ge 2.
\end{align*}

Fix $t\in[\varepsilon_{n+1},\, \varepsilon_n]$. Then
\begin{align*}
&\frac{1}{\eta(s)(s-1)2^{n+1}\ln^{s-1}(n+3)}\le2^{-n-1} r_{n+1}=\varepsilon_{n+1}\le t\\
&\le \varepsilon_n=2^{-n}r_n\le\frac{1}{\eta(s)(s-1)2^n\ln^{s-1}(n+1)}.
\end{align*}
Since $2^{n+1}\ln^{s-1}(n+3)$  and $2^{n-1}\ln^{s-1}(n+1)$ are  comparable with $2^n\ln^{s-1}n$ for large $n$ we can take positive constants $b_1,\, b_2$ such that
\begin{align}
&\frac{b_1}{2^n\ln^{s-1}n}\le t\le \frac{b_2}{2^n\ln^{s-1}n}.\label{efjkvefjobvenbo}
\end{align}
By Lemma \ref{decvefjvjjo} we have
\begin{align}
&r_{n+1}\le |K(t)|\le 4r_n,
\end{align}
and so there are two positive constants $d_1,\, d_2$ with
\begin{align*}
&\frac{d_1}{\ln^{s-1}n}\le |K(t)|\le \frac{d_2}{\ln^{s-1}n}.
\end{align*}
By \eqref{efjkvefjobvenbo} we obtain
\begin{align*}
&\ln b_1-n\ln2-(s-1)\ln \ln n\le \ln t\le \ln b_2-n\ln 2-(s-1)\ln \ln n,
\end{align*}
which gives for some constants $L_1,\, L_2$
\begin{align*}
&L_1 n\le n\ln2+(s-1)\ln \ln n\le \ln\frac{b_2}{t},\\
&\ln\frac{b_1}{t}\le n\ln2+(s-1)\ln \ln n\le L_2 n.
\end{align*}
So
\begin{align*}
&\frac{1}{L_2} \ln\frac{b_1}{t}\le n\le\frac{1}{L_1} \ln\frac{b_2}{t}.
\end{align*}
This implies
\begin{align*}
&\left(\ln \Big(\frac{1}{L_2} \ln\frac{b_1}{t}\Big)\right)^{s-1}\le \ln^{s-1} n\le\left(\ln \Big(\frac{1}{L_1} \ln\frac{b_2}{t}\Big)\right)^{s-1}.
\end{align*}
Finally we can find $ c_1,\, c_2$ and $b$ such that
\begin{align*}
&c_1 (\ln\ln(b/t))^{1-s}\le |K(t)|\le c_2  (\ln\ln(b/t))^{1-s}.
\end{align*}
\end{proof}

All Cantor sets have been constructed on an interval $[0,1]$ so far. But we can construct Cantor sets in $[0,1]^N$ as a cartesian product. But having $|K(t)|=\varphi(t)$ we have
$|K^N(t)|\le |K(t)|^N \le \varphi^N(t)$. In Lemmas \ref{djkojkovnfovn}, \ref{vjvjbvjvvjvjj} and \ref{dsjkvbdjovdjvj} we essentially get nothing new; the behavior of $\varphi(t)$
stays qualitatively the same.

\section{Examples of almost compact embeddings}

\begin{example}\label{kidvfuiovhfv}
Assume that $\Omega, p(\cdot),q(\cdot)$ satisfy \eqref{doicdocjdoij} and \eqref{eovjoievjvj}.
Let $K\subset\Omega$ be compact with zero measure and let $\varphi$ be given by \eqref{efovhevhiobh}. Suppose that $\varphi(t)\le Ct^s$ for some $C>0$ and $s\in (0,N]$. Assume that $\omega:(0,\diam(\Omega)]\rightarrow (0,\infty)$ satisfies
\begin{enumerate}[\rm(i)]
\item $\omega$ is decreasing;\label{dwjcjevebjlj}\\
\item $\lim\limits_{t\rightarrow 0_+}\omega(t)=\infty$;\label{djcdkvcdvdk}\\
\item $\lim\limits_{t\rightarrow 0_+}\frac{\ln(1/t)}{\omega(t)}=\infty$;\label{fklbgkpnjknk}\\
\item $\frac{1}{p(x)-q(x)}\le c\omega(\d_K(x))$ for come $c>0$.\label{efefkbngtkobhjri}
\end{enumerate}
Then $L^{p(\cdot)}(\Omega)\overset{\ast}{\hookrightarrow} L^{q(\cdot)}(\Omega)$.
\end{example}
\begin{proof}
Let $a>1$. Clearly, $s(x)\le c\omega(\d_K(x))$ by \eqref{efefkbngtkobhjri}. Using \eqref{fklbgkpnjknk} we have
\begin{align}
&\lim_{y\rightarrow\infty}\frac{\ln(1/\omega^{-1}(y))}{y}=\lim_{t\rightarrow 0_+}\frac{\ln(1/t)}{\omega(t)}=\infty.\label{vjdsovjovjodojvjvdo}
\end{align}
Then
\begin{align*}
&\int_{\omega_0}^\infty \varphi(\omega^{-1}(y))a^y dy\le C\int_{\omega_0}^\infty (\omega^{-1}(y))^sa^y dy=C\int_{\omega_0}^\infty \e^{s\ln(\omega^{-1}(y))}a^y dy\\
&= C\int_{\omega_0}^\infty\e^{-s\frac{\ln(1/{\omega^{-1}(y)})}{y}\ y+y\ln a} dy= C\int_{\omega_0}^\infty\e^{-y(s\frac{\ln(1/{\omega^{-1}(y)})}{y}-\ln a)} dy
=I.
\end{align*}
By \eqref{vjdsovjovjodojvjvdo} we have $\frac{\ln(1/\omega^{-1}(y))}{y}\rightarrow \infty$ for $y\rightarrow\infty$ and so, $s\frac{\ln(1/{\omega^{-1}(y)})}{y}-\ln a\ge 1$ for large $y$ which implies $I<\infty$. Now, Lemma \ref{fkvjkkvjovji} gives
$L^{p(\cdot)}(\Omega)\overset{\ast}{\hookrightarrow} L^{q(\cdot)}(\Omega)$.
\end{proof}

\begin{example}\label{edjceddkvcjrvjrvj}
Assume that $\Omega, p(\cdot),q(\cdot)$ satisfy \eqref{doicdocjdoij} and \eqref{eovjoievjvj}.
Let $K\subset\Omega$ be compact with zero measure and let $\varphi$ be given by \eqref{efovhevhiobh} and $\varphi(t)\le C(\ln (B/t))^{1-s}$ for some $B,C>0$, $s>1$ and $t\in(0,|\Omega|)$. Assume that $\omega:(0,\diam(\Omega)]\rightarrow (0,\infty)$ satisfies
\begin{enumerate}[\rm(i)]
\item $\omega$ is decreasing;\label{giopnjgopnjoph}\\
\item $\lim\limits_{t\rightarrow 0_+}\omega(t)=\infty$;\label{dkvjdvjdvdooo}\\
\item $\lim\limits_{t\rightarrow 0_+}\frac{\ln\ln(1/t)}{\omega(t)}=\infty$;\label{opghgkonpjkopn}\\
\item $\frac{1}{p(x)-q(x)}\le c \omega(\d_K(x))$ for some $c>0$.\label{opgyopnjnjo}
\end{enumerate}
Then $L^{p(\cdot)}(\Omega)\overset{\ast}{\hookrightarrow} L^{q(\cdot)}(\Omega)$.
\end{example}
\begin{proof}
Let $a>1$. Clearly, $s(x)\le c\omega(\d_K(x))$ by \eqref{opgyopnjnjo}. Using \eqref{opghgkonpjkopn} we have
\begin{align}
&\lim_{y\rightarrow\infty}\frac{\ln\ln(1/\omega^{-1}(y))}{y}=\lim_{t\rightarrow 0_+}\frac{\ln\ln(1/t)}{\omega(t)}=\infty.\label{sklccjdwvjdvjv}
\end{align}
Then
\begin{align*}
&\int_{\omega_0}^\infty \varphi(\omega^{-1}(y))a^y dy\le C\int_{\omega_0}^\infty (\ln(B/\omega^{-1}(y))^{1-s}a^y dy\\
&=C\int_{\omega_0}^\infty \e^{(1-s)\ln\ln(B/\omega^{-1}(y))}a^y dy=C\int_{\omega_0}^\infty \e^{(1-s)\frac{\ln\ln(B/\omega^{-1}(y))}{y}\ y+y\ln a} dy\\
&=C\int_{\omega_0}^\infty \e^{y((1-s)\frac{\ln\ln(B/\omega^{-1}(y))}{y}+\ln a)} dy =I.
\end{align*}
By \eqref{sklccjdwvjdvjv} we have $\frac{\ln\ln(1/\omega^{-1}(y))}{y}\rightarrow \infty$ for $y\rightarrow\infty$. Moreover, $1-s<0$ and so, $(1-s)\frac{\ln\ln(B/\omega^{-1}(y))}{y}+\ln a\le -1$ for large $y$ which implies $I<\infty$. Now, Lemma \ref{fkvjkkvjovji} gives
$L^{p(\cdot)}(\Omega)\overset{\ast}{\hookrightarrow} L^{q(\cdot)}(\Omega)$.
\end{proof}

\begin{example}\label{efkvjpfevfekjvfekvkj}
Assume that $\Omega, p(\cdot),q(\cdot)$ satisfy \eqref{doicdocjdoij} and \eqref{eovjoievjvj}.
Let $K\subset\mathbb{R}^N$ be  compact with zero measure, $\varphi$ be given by \eqref{efovhevhiobh} and $\varphi(t)\le C(\ln\ln (B/t))^{1-s}$ for some $B,C>0$, $s>1$ and $t\in(0,|\Omega|)$. Assume that $\omega:(0,\diam(\Omega)]\rightarrow (0,\infty)$ satisfies
\begin{enumerate}[\rm(i)]
\item $\omega$ is decreasing;\label{hjiophjhjjopjop}\\
\item $\lim\limits_{t\rightarrow 0_+}\omega(t)=\infty$;\label{djvdoivivuiueio}\\
\item $\lim\limits_{t\rightarrow 0_+}\frac{\ln\ln\ln(1/t)}{\omega(t)}=\infty$;\label{cvpbgfpmgnmjnmlp}\\
\item $\frac{1}{p(x)-q(x)}\le c\omega(\d_K(x))$ for some $c>0$.\label{kldfbbkkgnknk}
\end{enumerate}
Then $L^{p(\cdot)}(\Omega)\overset{\ast}{\hookrightarrow} L^{q(\cdot)}(\Omega)$.
\end{example}
\begin{proof}
Let $a>1$. Clearly, $s(x)\le c\omega(\d_K(x))$ by \eqref{kldfbbkkgnknk}. Using \eqref{cvpbgfpmgnmjnmlp} we have
\begin{align}
&\lim_{y\rightarrow\infty}\frac{\ln\ln\ln(1/\omega^{-1}(y))}{y}=\lim_{t\rightarrow 0_+}\frac{\ln\ln\ln(1/t)}{\omega(t)}=\infty.\label{dlvkjdvjdvjdvj}
\end{align}
Then
\begin{align*}
&\int_{\omega_0}^\infty \varphi(\omega^{-1}(y))a^y dy\le C\int_{\omega_0}^\infty (\ln\ln(B/\omega^{-1}(y))^{1-s}a^y dy\\
&=C\int_{\omega_0}^\infty \e^{(1-s)\ln\ln\ln(B/\omega^{-1}(y))}a^y dy=C\int_{\omega_0}^\infty \e^{(1-s)\frac{\ln\ln\ln(B/\omega^{-1}(y))}{y}\ y+y\ln a} dy\\
&=C\int_{\omega_0}^\infty \e^{y((1-s)\frac{\ln\ln\ln(B/\omega^{-1}(y))}{y}+\ln a)} dy =I.
\end{align*}
By \eqref{dlvkjdvjdvjdvj} we have $\frac{\ln\ln\ln(1/\omega^{-1}(y))}{y}\rightarrow \infty$ for $y\rightarrow\infty$. Moreover, $1-s<0$ and so, $(1-s)\frac{\ln\ln\ln(B/\omega^{-1}(y))}{y}+\ln a\le -1$ for large $y$ which implies $I<\infty$. Now, Lemma~\ref{fkvjkkvjovji} gives
$L^{p(\cdot)}(\Omega)\overset{\ast}{\hookrightarrow} L^{q(\cdot)}(\Omega)$.
\end{proof}

\section{Compact embeddings between variable Sobolev and variable  Lebesgue spaces}

First of all we establish a necessary condition for an embedding to be compact.

\begin{lemma}\label{dkvcefjvhnj}
Let $B_r=B(0,r)$ denote the ball in $\mathbb{R}^N$ centered at $0$ with radius $r$. Assume $M\subset B_r$ and $s\in(0,r)$ are such that $|B_r\setminus B_s|\le |M|$. Suppose that $\varphi:(0,r]\rightarrow\mathbb{R}$ is non-negative and non-increasing and set $\psi(x)=\varphi(|x|)$, $x\in B_r$. Then
\begin{align*}
&\int_M \psi(x)dx\ge \int_{B_r\setminus B_s} \psi(x)dx.
\end{align*}
\end{lemma}
\begin{proof}
By the assumption  $|B_r\setminus B_s|\le |M|$ we have
\begin{align*}|(B_r\setminus B_s)\setminus M|+|(B_r\setminus B_s)\cap M|&=|(B_r\setminus B_s)|\le |M|\\
&=|(B_r\setminus B_s)\cap M|+|M\cap B_s|.
\end{align*}
and consequently
\[|M\cap B_s| \ge |(B_r\setminus B_s)\setminus M|.
\]
By the assumptions on  $\psi$ we have
$\psi(x)\ge \psi(y)$ for every $x\in B_s$ and every $y\in  B_r\setminus B_s$.
This implies
\begin{align*}
\int_{M}\psi(x)dx&=\int_{M\cap B_s}\psi(x)dx+\int_{(B_r\setminus B_s)\cap M}\psi(x)dx\\
&\ge \frac{|M\cap B_s|}{|(B_r\setminus B_s)\setminus M|}\int_{(B_r\setminus B_s)\setminus M}\psi(x)dx +\int_{(B_r\setminus B_s)\cap M}\psi(x)dx\\
&\ge \int_{(B_r\setminus B_s)\setminus M}\psi(x)dx +\int_{(B_r\setminus B_s)\cap M}\psi(x)dx\\
&=\int_{B_r\setminus B_s}\psi(x)dx,
\end{align*}
which finishes the proof.
\end{proof}

\begin{theorem}\label{dklvnfkobgb}
Let $p,\, q\in\mathcal{E}(\Omega)$, $1\le p(x)\le p_{+}<N$ on $\Omega$ and let $p(\cdot)$ satisfy \eqref{devjfjkovf}, $1\le q(x)\le p^\#(x)$ and let $M=\{x\in \Omega;\,\, p^\#(x)=q(x)\}$. Assume
\begin{align*}
&W^{1,p(\cdot)}(\Omega)\hookrightarrow\hookrightarrow L^{q(\cdot)}(\Omega).
\end{align*}
Then $|M|=0$.
\end{theorem}
\begin{proof}
Suppose $|M|>0$. Let $x_0\in \Omega$ be a point of Lebesgue density of $M$. Given $\varepsilon>0$ denote
\begin{align*}
&p_-^\varepsilon=\inf\{p(x);x\in B(x_0,\varepsilon)\}, \,\,\,\, p_{+}^\varepsilon=\sup\{p(x);x\in B(x_0,\varepsilon)\},
\end{align*}
and define a function $u_\varepsilon$ by
\begin{align*}
&u_\varepsilon(x)=\varepsilon^{\frac{p_-^\varepsilon-N}{p_-^\varepsilon}}(1-{|x-x_0|}/{\varepsilon})\chi_{B(x_0,\varepsilon)}(x).
\end{align*}
Clearly,
\begin{align*}
&|\nabla u_\varepsilon(x)|=\varepsilon^{\frac{p_-^\varepsilon-N}{p_-^\varepsilon}}\ {1}/{\varepsilon}\ \chi_{B(x_0,\varepsilon)}(x)=\varepsilon^{-\frac{N}{p_-^\varepsilon}}\ \chi_{B(x_0,\varepsilon)}(x).
\end{align*}

First we prove that the set $\{u_\varepsilon\}$ is bounded in $W^{1,p(\cdot)}(\Omega)$ for $\varepsilon\le 1$. Plainly,
\begin{align*}
&\int_\Omega |u_\varepsilon(x)|^{p(x)}dx=\int_{B(x_0,\varepsilon)}  \varepsilon^{\frac{p_-^\varepsilon-N}{p_-^\varepsilon}p(x)}(1-{|x-x_0|}/{\varepsilon})^{p(x)}dx\\
&\le \int_{B(x_0,\varepsilon)}  \varepsilon^{p(x)}\varepsilon^{-\frac{N}{p_-^\varepsilon}p(x)}dx\le  \int_{B(x_0,\varepsilon)}\varepsilon^{-\frac{N}{p_-^\varepsilon}p(x)}dx:=I_\varepsilon.
\end{align*}
Moreover,
\begin{align*}
&\int_\Omega |\nabla u_\varepsilon(x)|^{p(x)}dx\le  \int_{B(x_0,\varepsilon)}\varepsilon^{-\frac{N}{p_-^\varepsilon}p(x)}dx=:I_\varepsilon.
\end{align*}
Now,
\begin{align*}
&I_\varepsilon=\int_{B(x_0,\varepsilon)}\varepsilon^{-\frac{N}{p_-^\varepsilon}(p(x)-p_-^\varepsilon)}\ \varepsilon^{-N}dx=
\int_{B(x_0,\varepsilon)}\e^{-\frac{N}{p_-^\varepsilon}(p(x)-p_-^\varepsilon)\ln \varepsilon}\ \varepsilon^{-N}dx\\
&=\int_{B(x_0,\varepsilon)}\e^{\frac{N}{p_-^\varepsilon}(p(x)-p_-^\varepsilon)\ln (1/\varepsilon)}\ \varepsilon^{-N}dx.
\end{align*}
From the log-Lipschitz condition \eqref{devjfjkovf} we have
\begin{align*}
(p(x)-p_-^\varepsilon)\ln (1/\varepsilon)\le C
\end{align*}
and so
\begin{align*}
&I_\varepsilon\le \int_{B(x_0,\varepsilon)}\e^{\frac{CN}{p_-^\varepsilon}}\varepsilon^{-N}dx\le \e^{CN}\int_{B(x_0,\varepsilon)}\varepsilon^{-N}dx=:A.
\end{align*}
This immediately implies that $\|u_\varepsilon\|_{W^{1,p(\cdot)}(\Omega)}$ is bounded  for $\varepsilon\le 1$.

Now fix $\varepsilon_0$ such that for all $\varepsilon\le \varepsilon_0$ we have
\begin{align*}
|B(x_0,\varepsilon)\cap M|\ge |B(x_0,7\varepsilon/8)|.
\end{align*}

Fix for a moment $\varepsilon\le \varepsilon_0$. Then
\begin{align*}
&A_\varepsilon:=\int_{(B(x_0,\varepsilon)\setminus B(x_0,(3/4)\varepsilon))\cap M} (1-|x-x_0|/\varepsilon)^{p^\#(x)}\Big(\varepsilon^{\frac{p_-^\varepsilon-N}{p_-^\varepsilon}}\Big)^{p^\#(x)}dx\\
&\ge \int_{(B(x_0,\varepsilon)\setminus B(x_0,(3/4)\varepsilon))\cap M} (1-|x-x_0|/\varepsilon)^{p^\#(x)}\varepsilon^{\frac{p_-^\varepsilon-N}{p_-^\varepsilon}(p_-^\varepsilon)^\#}dx\\
&=\int_{(B(x_0,\varepsilon)\setminus B(x_0,(3/4)\varepsilon))\cap M} (1-|x-x_0|/\varepsilon)^{p^\#(x)}\varepsilon^{\frac{p_-^\varepsilon-N}{p_-^\varepsilon}\frac{Np_-^\varepsilon}{N-p_-^\varepsilon}}dx\\
&=\int_{(B(x_0,\varepsilon)\setminus B(x_0,(3/4)\varepsilon))\cap M} (1-|x-x_0|/\varepsilon)^{p^\#(x)}\varepsilon^{-N}dx\\
&\ge \int_{(B(x_0,\varepsilon)\setminus B(x_0,(3/4)\varepsilon))\cap M} (1-|x-x_0|/\varepsilon)^{(p_{+}^\varepsilon)^\#}\varepsilon^{-N}dx=:B_\varepsilon.
\end{align*}
By Lemma \ref{dkvcefjvhnj} we have
\begin{align*}
&B_\varepsilon\ge \int_{(B(x_0,\varepsilon)\setminus B(x_0,(7/8)\varepsilon)} (1-|x-x_0|/\varepsilon)^{(p_+^\varepsilon)^\#}\varepsilon^{-N}dx\\
&=\sigma_N\varepsilon^{-N}\int_{(7/8)\varepsilon}^{\varepsilon} (1-r/\varepsilon)^{(p_+^\varepsilon)^\#}r^{N-1}dr\\
&\ge \sigma_N\varepsilon^{-N}\int_{(14/16)\varepsilon}^{(15/16)\varepsilon} (1-(15\varepsilon/16)/\varepsilon)^{(p_+^\varepsilon)^\#}r^{N-1}dr\\
&=\sigma_N  (1/16)^{(p_+^\varepsilon)^\#}\   \varepsilon^{-N}\int_{(14/16)\varepsilon}^{(15/16)\varepsilon}r^{N-1}dr=:K.
\end{align*}
where $\sigma_N$ denotes the area of  the $N$-dimensional unit sphere.

Denote $\varepsilon_n=(3/4)^n \varepsilon_0$ and consider the corresponding sequence $u_{\varepsilon_n}(x)$. Let $m>n$. Then $u_m(x)=0$ for $x\in B(x_0,\varepsilon_n)\setminus B(x_0,\varepsilon_m)$ and so,
\begin{align*}
&\int_\Omega |u_m(x)-u_n(x)|^{q(x)}dx=\int_{B(x_0,\varepsilon_n)} |u_m(x)-u_n(x)|^{q(x)}dx\\
&\ge \int_{B(x_0,\varepsilon_n)\setminus B(x_0,\varepsilon_m)} |u_n(x)|^{q(x)}dx\ge\int_{(B(x_0,\varepsilon_n)\setminus B(x_0,(3/4)\varepsilon_n))\cap M} |u_n(x)|^{q(x)} dx\\
&\ge \int_{(B(x_0,\varepsilon_n)\setminus B(x_0,(3/4)\varepsilon_n))\cap M} |u_n(x)|^{p^\#(x)} dx\ge K.
\end{align*}
Hence, there is a constant $L>0$ such that
\begin{align*}
&\|u_m-u_n\|_{L^{q(\cdot)}(\Omega)}\ge L
\end{align*}
and the embedding $W^{1,p(\cdot)}(\Omega)\hookrightarrow\hookrightarrow L^{q(\cdot)}(\Omega)$ is not compact.
\end{proof}

The next lemma is proved in \cite{DHHR} (see Corollary 8.3.2.).

\begin{lemma}\label{jnskljvnsljkv}
Let $\Omega\in\mathcal{C}^{0,1}$, $p, \, q\in\mathcal{E}(\Omega)$ and suppose $p(\cdot)$ satisfies the $log$-H\"older condition \eqref{devjfjkovf}. Assume that for all $x\in \Omega$
\begin{align*}
&1\le p(x)\le p_{+}<N.
\end{align*}
Then $W^{1,p(\cdot)}(\Omega)\hookrightarrow L^{p^\#(\cdot)}(\Omega)$ where $p^\#(x)$ is given in \eqref{jvhbjovhobvh}.
\end{lemma}

\begin{theorem}\label{oiwdvhovwoiio}
Let $\Omega\in\mathcal{C}^{0,1}$, $p,\, q\in\mathcal{E}(\Omega)$ and let $p(\cdot)$ satisfy the $log$-H\"older condition \eqref{devjfjkovf}.  Assume that for all $x\in \Omega$,
\begin{align*}
&1\le p(x)\le p_{+}<N,\ \ 1\le q(x)\le p^\#(x)
\end{align*}
where $p^\#(x)$ is given in \eqref{jvhbjovhobvh}. Let $K\subset \overline{\Omega}$ be compact, $|K|=0$ and denote $\varphi(t)=|K(t)|$.
Let $\omega:(0,\diam(\Omega)]\rightarrow\mathbb{R}$ be a decreasing continuous non-negative function, $\omega_0:=\omega(\diam(\Omega))$. Suppose that $\omega(\cdot)$ satisfies
\begin{align*}
&\frac{1}{p^\#(x)-q(x)}\le c\ \omega(\d_K(x)), \ x\in \Omega,\\
&\int_{\omega_0}^{\omega(0)} \varphi(\omega^{-1}(y))a^y dy<\infty \ \text{for all}\ a>1.
\end{align*}
Then $W^{1,p(\cdot)}(\Omega)\hookrightarrow\hookrightarrow L^{q(\cdot)}(\Omega)$.
\end{theorem}
\begin{proof}
Lemmas \ref{jnskljvnsljkv} and \ref{fkvjkkvjovji} give
\begin{align*}
&W^{1,p(\cdot)}(\Omega)\hookrightarrow L^{p^{\#}}(\Omega)\overset{\ast}{\hookrightarrow} L^{q(\cdot)}(\Omega).
\end{align*}
Now Proposition \ref{wiuhvuvhbhobi} finishes the proof.
\end{proof}

As an application we introduce the following several examples. The first one is in fact proved in \cite{MOSS} (see Theorem 3.4) but we obtain it as an easy consequence of the previous theorem.
Let $\varphi$ and $q$ denote in following three examples the same as in Theorem \ref{oiwdvhovwoiio}.
\begin{example}
Let $\Omega\in\mathcal{C}^{0,1}$ and $p(\cdot):\Omega\rightarrow\mathbb{R}$ satisfy \eqref{devjfjkovf}.  Assume
\begin{align*}
&1\le p(x)\le p_{+}<N.
\end{align*}
Let $K\subset\overline{\Omega}$ be a compact set with zero measure and suppose that $\varphi(t)\le Ct^s$ for some $C>0$ and $s\in (0,N]$. Assume that $\omega:(0,\diam(\Omega)]\rightarrow (0,\infty)$ satisfies
\begin{enumerate}[\rm(i)]
\item $\omega$ is decreasing and continuous;\\
\item $\lim\limits_{t\rightarrow 0_+}\omega(t)=\infty$;\\
\item $\lim\limits_{t\rightarrow 0_+}\frac{\ln(1/t)}{\omega(t)}=\infty$;\\
\item $\frac{1}{p^{\#}(x)-q(x)}\le c\omega(\d_K(x))$ for come $c>0$.
\end{enumerate}
Then $W^{1,p(\cdot)}(\Omega){\hookrightarrow\hookrightarrow} L^{q(\cdot)}(\Omega)$.
\end{example}
\begin{proof}
It suffices to use Theorem \ref{oiwdvhovwoiio} and Example \ref{kidvfuiovhfv}.
\end{proof}

\begin{example}
Let $\Omega\in\mathcal{C}^{0,1}$ and $p(\cdot):\Omega\rightarrow\mathbb{R}$ satisfy \eqref{devjfjkovf}.  Assume
\begin{align*}
&1\le p(x)\le p_{+}<N.
\end{align*}
Let $K\subset\overline{\Omega}$ be a compact set with zero measure and suppose that $\varphi(t)\le C(\ln (B/t))^{1-s}$ for some $B,C>0$, $s>1$ and $t\in(0,|\Omega|)$. Assume that $\omega:(0,\diam(\Omega)]\rightarrow (0,\infty)$ satisfies
\begin{enumerate}[\rm(i)]
\item $\omega$ is decreasing and continuous;\\
\item $\lim\limits_{t\rightarrow 0_+}\omega(t)=\infty$;\\
\item $\lim\limits_{t\rightarrow 0_+}\frac{\ln\ln(1/t)}{\omega(t)}=\infty$;\\
\item $\frac{1}{p^{\#}(x)-q(x)}\le c \omega(\d_K(x))$ for some $c>0$.
\end{enumerate}
Then $W^{1,p(\cdot)}(\Omega){\hookrightarrow\hookrightarrow} L^{q(\cdot)}(\Omega)$.
\end{example}
\begin{proof}
Use Theorem \ref{oiwdvhovwoiio} and Example \ref{edjceddkvcjrvjrvj}.
\end{proof}

\begin{example}
Let $\Omega\in\mathcal{C}^{0,1}$ and $p(\cdot):\Omega\rightarrow\mathbb{R}$ satisfy \eqref{devjfjkovf}.  Assume
\begin{align*}
&1\le p(x)\le p_{+}<N.
\end{align*}
Let $K\subset\overline{\Omega}$ be  of zero measure and $\varphi(t)\le C(\ln\ln (B/t))^{1-s}$ for some $B,C>0$, $s>1$ and $t\in(0,|\Omega|)$. Assume that $\omega:(0,\diam(\Omega)]\rightarrow (0,\infty)$ satisfies
\begin{enumerate}[\rm(i)]
\item $\omega$ is decreasing and continuous;\\
\item $\lim\limits_{t\rightarrow 0_+}\omega(t)=\infty$;\\
\item $\lim\limits_{t\rightarrow 0_+}\frac{\ln\ln\ln(1/t)}{\omega(t)}=\infty$;\\
\item $\frac{1}{p^{\#}(x)-q(x)}\le c \omega(\d_K(x))$ for some $c>0$.
\end{enumerate}
Then $W^{1,p(\cdot)}(\Omega){\hookrightarrow\hookrightarrow} L^{q(\cdot)}(\Omega)$.
\end{example}
\begin{proof}
It suffices to use Theorem \ref{oiwdvhovwoiio} and Example \ref{efkvjpfevfekjvfekvkj}.
\end{proof}

To conclude we remark that the construction of Cantor sets could be refined by adding some more logarithms to give additional examples.

\

\subsection*{Acknowledgements}
The second and the third author of this research were supported by the grant P201-18-00580S of the Grant Agency of the Czech Republic. The second author has been partially supported by Shota Rustaveli National Science Foundation of Georgia (SRNSFG) [grant number FR17-589] and RVO:67985840.

We would like to thank anonymous referees for their critical reading of the paper and for numerous useful comments and suggestions.


\end{document}